\documentclass{elsarticle}
\usepackage{jcomp}
\usepackage{lineno}
\usepackage{soul}

\usepackage[utf8]{inputenc}
\usepackage{hyperref}
\usepackage{mathtools}
\usepackage{amsthm}
\usepackage{tikz}
\usepackage{cleveref}
\usepackage{graphicx}
\usepackage{booktabs}
\usepackage{subcaption}
\usepackage{overpic} % To add text on Figures

\DeclareMathAccent{\svec}{\mathord}{letters}{126}
\DeclareMathOperator{\diag}{diag}
\newcommand{\mat}[1]{\ensuremath{\mathsf{#1}}}
\newcommand{\vect}[1]{\ensuremath{\boldsymbol{\mathbf{#1}}}}
\newcommand{\avg}[1]{\ensuremath{\left\{\hspace*{-3pt}\left\{#1\right\}\hspace*{-3pt}\right\}}}
\newcommand{\SC}{\ensuremath{\text{S}}}
\newcommand{\RS}{\ensuremath{\star}}

\newenvironment{lequation}{%
    \begin{linenomath*}\begin{equation}%
}{%
    \end{equation}\end{linenomath*}%
}
\newenvironment{lequation*}{%
    \begin{linenomath*}\begin{equation*}%
}{%
    \end{equation*}\end{linenomath*}%
}

\newtheorem{theorem}{Theorem}

\newcommand{\one}[1]{#1}
\newcommand{\two}[1]{#1}
\newcommand{\three}[1]{#1}

\renewcommand{\st}[1]{}

\begin{document}
%---------------------------------------------------------------------------------------------------

%---------------------------------------------------------------------------------------------------
\begin{frontmatter}
%---------------------------------------------------------------------------------------------------

\title{A flux-differencing formulation with Gauss nodes}

\journal{Journal of Computational Physics}

\author[1,2]{Andr\'es \snm{Mateo-Gab\'in}\corref{cor1}}
\author[3]{Andr\'es M. \snm{Rueda-Ram\'irez}}
\author[1,2]{Eusebio \snm{Valero}}
\author[1,2]{Gonzalo \snm{Rubio}}

\cortext[cor1]{Corresponding author: andres.mgabin@upm.es}
\address[1]{Universidad Polit\'ecnica de Madrid, School of Aeronautics, Madrid, Spain}
\address[2]{Universidad Polit\'ecnica de Madrid, Center for Computational Simulation, Madrid, Spain}
\address[3]{University of Cologne, Department of Mathematics and Computer Science, Cologne, Germany}

%\begin{keyword}
%    discontinuous Galerkin \sep high order \sep entropy stability \sep shock capturing \sep summation-by-parts operators \sep flux differencing
%\end{keyword}

\end{frontmatter}
% \linenumbers

%---------------------------------------------------------------------------------------------------
\section{Introduction}
%---------------------------------------------------------------------------------------------------

Among the different spatial discretization frameworks, the Discontinuous Galerkin Spectral Element Method (DGSEM)~\cite{black1999conservative,Kopriva2009} approximates a function as a set of discontinuous, high-order polynomials defined in a tessellation of the spatial domain. 
The DGSEM is a collocation method that represents the solution as a combination of Lagrange polynomials that shares the nodes with a certain quadrature rule, usually Gauss or Gauss-Lobatto.

Fisher \three{et al.} showed in~\cite{Fisher2012} that diagonal-norm summation-by-parts (SBP) derivative operators can be rewritten in telescopic form. In this case, the action of the operator resembles a finite volume scheme, ensuring local conservation in the sense of Lax-Wendroff. In the context of the DGSEM, the derivative operator with Gauss-Lobatto nodes falls in this category, whereas the use of Gauss nodes leads to a generalized SBP operator~\one{\cite{DelRey2014,hicken2016multidimensional,del2019extension,Chan2018,chan2019efficient}}.

In this work we show that the generalized formulation of Chan~\cite{chan2019efficient} also admits a flux-differencing form and thus, from a more practical standpoint, entropy-stable schemes can be generated by introducing certain dissipative numerical fluxes at the interfaces of the sub-elements~\three{\cite{Fisher2013,Hennemann2021,RUEDARAMIREZ2022subcell}}.

%---------------------------------------------------------------------------------------------------
\section{High-order DGSEM on Gauss nodes} \label{sec:theory}
%---------------------------------------------------------------------------------------------------

\begin{figure}
\centering
\begin{tikzpicture}[thick, scale=0.6]
    \draw (-4.4,0) -- (4.4,0);
    \draw (-4.0,0.5em) -- (-4.0,-0.5em) node[above=0.6em] {$\bar{\xi}_0$};
    \draw (-2.6086,0.5em) -- (-2.6086,-0.5em) node[above=0.6em] {$\bar{\xi}_1$};
    \draw (0.0,0.5em) -- (0.0,-0.5em) node[above=0.6em] {$\bar{\xi}_2$};
    \draw (2.6086,0.5em) -- (2.6086,-0.5em) node[above=0.6em] {$\bar{\xi}_3$};
    \draw (4.0,0.5em) -- (4.0,-0.5em) node[above=0.6em] {$\bar{\xi}_4$};
    \filldraw (-3.4445,0.0) circle (3pt) node[below] {$\xi_0$};
    \filldraw (-1.3599,0.0) circle (3pt) node[below] {$\xi_1$};
    \filldraw (1.3599,0.0) circle (3pt) node[below] {$\xi_2$};
    \filldraw (3.4445,0.0) circle (3pt) node[below] {$\xi_3$};
\end{tikzpicture}
\caption{Main,~$\xi_i$, and complementary,~$\bar{\xi}_i$, grids for~$N=3$.}
\label{fig:grid}
\end{figure}
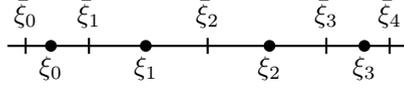

Beginning with a one-dimensional grid where~$\xi_i$ and~$\omega_i$ are the nodes and weights of the Gauss quadrature rule of order~$N$ ($i \in [0,N]$) in the reference domain,~$E$, we define a complementary grid~$\bar{\xi}_i$ ($i \in [0,N+1]$) as in \cref{fig:grid},
\begin{lequation*}
    \bar{\xi}_0 = -1, \quad \bar{\xi}_i = \bar{\xi}_{i-1} + \omega_i, \quad \bar{\xi}_{N+1} = 1.
\end{lequation*}

Fisher \three{et al.} showed in~\cite{Fisher2012} that SBP derivative operators,~$\mat{D}$, can be rewritten in telescopic form. In this case, the action of the operator resembles a finite volume scheme,
\begin{lequation} \label{eq:fv}
    \left(\partial_\xi f\right)_i \approx \sum_{k=0}^N D_{ik} f_k = \frac{\bar{f}_{i+1} - \bar{f}_i}{\omega_i}.
\end{lequation}

As stated in the introduction, the derivative operator of the DGSEM with Gauss-Lobatto nodes falls in this category, while we can use the developments of Chan~\cite{Chan2018,chan2019efficient} if Gauss nodes are used. Following Chan's notation, the DGSEM discretization of a generic one-dimensional conservation law, $u_t + f_x = 0$, is
\begin{lequation} \label{eq:split_chan}
    \renewcommand{\arraystretch}{1.3}
    \mat{M}\two{\vect{u_t}} + \left[\mat{I} \quad \mat{V_f}^T\mat{B}\right]\left(2\mat{Q}\circ\mat{F^{\SC}}\right)\cdot\vect{1} + \mat{V_f}^T\mat{B}\left[\begin{array}{c}
        f^{\RS}_L - f_{\tilde{L}} \\
        f^{\RS}_R - f_{\tilde{R}}
    \end{array}\right] = 0.
\end{lequation}
The mass matrix,~$\mat{M}=\diag(J\omega_i)$, contains the Jacobians and weights of the quadrature rule in each mesh element. \one{The projection of the interior values to the faces is represented by~$\mat{V_f}$,} and~$\mat{B}=\diag(-1,1)$. \one{The matrix}~$\mat{Q}$ is a generalized SBP operator representing the integral of the derivative of the basis functions,
\one{\begin{lequation*}
    2\mat{Q} = \left[\begin{array}{cc}
        \mat{S} & \mat{V_f}^T\mat{B} \\
        -\mat{B}\mat{V_f} & \mat{B} \\
    \end{array}\right], \quad \mat{S} = 2\mat{M}\mat{D} - \mat{V_f}^T\mat{B}\mat{V_f}, \quad D_{ij} = l_j'(\xi_i), \quad \mat{V_f} = \left[\begin{array}{cccc}
        l_0(-1) & l_1(-1) & \cdots & l_N(-1) \\
        l_0(+1) & l_1(+1) & \cdots & l_N(+1) \\
    \end{array}\right],
\end{lequation*}}
\one{where $l_i(x)$ is the Lagrange polynomial associated to the node $\xi_i$. The numerical fluxes are indicated with a star ($f^{\RS}_L$ and~$f^{\RS}_R$), and} sub-indices with a tilde mean that the magnitude has been computed on so-called entropy projected variables,
\one{\begin{lequation*}
    f_{\tilde{L}} = f\left(u(v_L)\right), \quad f_{\tilde{R}} = f\left(u(v_R)\right), \quad \mat{V_f}\cdot\vect{v}(\vect{u}) = \left[\begin{array}{c}
        v_L \\ v_R \\
    \end{array}\right].
\end{lequation*}}
Finally, $\mat{F^{\SC}}$ is the matrix with all the combinations of entropy-conservative fluxes (including the entropy-projected values at the interfaces),~$F^{\SC}_{ij}=f^{\SC}(u_i, u_j)$, and~\two{$\mat{A} \circ \mat{B}$} is the Hadamard product \two{of matrices~$\mat{A}$ and~$\mat{B}$}. \one{More details} can be found in~\cite{Chan2018}.
 
Imposing the equality of~\cref{eq:split_chan} with a finite volume scheme,
\begin{lequation} \label{eq:split_chan_fv}
    \mat{M}\two{\vect{u_t}} + \mat{\Delta}\vect{\bar{f}} = 0,
\end{lequation}
where~$\mat{\Delta}\vect{\bar{f}}$ represents the flux-differencing formulation of \cref{eq:fv}, we find an expression for the different subcell fluxes,~$\vect{\bar{f}}$,
\begin{lequation} \label{eq:ssfv}
    \renewcommand{\arraystretch}{1.3}
    \mat{\Delta}\vect{\bar{f}} = \left[\mat{I} \quad \mat{V_f}^T\mat{B}\right]\left(2\mat{Q}\circ\mat{F^{\SC}}\right)\cdot\vect{1} + \mat{V_f}^T\mat{B}\left[\begin{array}{c}
        f^{\RS}_L - f_{\tilde{L}} \\
        f^{\RS}_R - f_{\tilde{R}}
    \end{array}\right].
\end{lequation}
With the additional constraints $\bar{f}_0 = f^{\RS}_L$ and $\bar{f}_{N+1} = f^{\RS}_R$, we can compute the value of~$\bar{f}_{i+1}$ from~$\bar{f}_i$ and \cref{eq:ssfv}, obtaining an expression \two{for the interface fluxes that resembles the one developed by Rueda-Ram\'irez et al.~\cite{Rueda2023} for the semi-discretization},
\begin{lequation} \label{eq:fc_definition}
\begin{gathered}
    \bar{f}_0 = f^{\RS}_L, \quad \bar{f}_{N+1} = f^{\RS}_R, \\
    \bar{f}_{i+1} = \bar{f}_i + \sum_{k=0}^N \one{S_{ik}}f^{\SC}_{ik} - l_i(-1)\left[f^{\SC}_{i\tilde{L}} - \sum_{k=0}^N l_k(-1)f^{\SC}_{\tilde{L}k} + f^{\RS}_L\right] + l_i(+1)\left[f^{\SC}_{i\tilde{R}} - \sum_{k=0}^N l_k(+1)f^{\SC}_{\tilde{R}k} + f^{\RS}_R\right], \quad i=0, \ldots, N,
\end{gathered}
\end{lequation}
\one{\st{where~$\mat{\hat{S}}= 2\mat{Q}-\mat{B}=\mat{Q}-\mat{Q}^T$ is a new skew-symmetric derivative matrix that includes the inner boundary contributions. }}We remark that \cref{eq:split_chan} is recovered when $\bar{f}_0$ and $\bar{f}_{i+1}$ are introduced in \cref{eq:split_chan_fv} by construction. \three{The last two terms of \cref{eq:fc_definition} represent a new addition with respect to eq. (3.9) of~\cite{Fisher2013}, and they couple the volume and surface integrals. These new operations result in a higher computational cost. However, the use of Gauss nodes also entails a higher accuracy that can compensate it as shown in~\cite{Rueda2023}.}

There is, however, a problem with~\cref{eq:fc_definition}. For a set of~$N+1$ Gauss nodes there are~$N+2$ complementary staggered fluxes, but we have~$N+3$ equations. The flux at the right face,~$\bar{f}_{N+1}$, can be computed from~$\bar{f}_{i+1}(i=N)$, but it also must be equal to~$f^{\RS}_R$. We can overcome this issue by proving that both equations are equivalent and our system is not over-constrained.
\two{\begin{theorem}
The set of equations~\eqref{eq:fc_definition} uniquely defines the interface fluxes of the subcell grid.
\end{theorem}}
\begin{proof}
Since the complementary flux~$\bar{f}_{i+1}$ is defined in terms of~$\bar{f}_i$, it is possible to define all the complementary fluxes in terms of the first one,~$\bar{f}_0=f^{\RS}_L$,
\begin{lequation} \label{eq:fc_definition_L}
    \bar{f}_{i+1} = f^{\RS}_L + \sum_{\alpha=0}^i\sum_{k=0}^N \one{S_{\alpha k}}f^{\SC}_{\alpha k} - \sum_{\alpha=0}^i l_{\alpha}(-1)\left(f^{\SC}_{\alpha\tilde{L}} - \sum_{k=0}^N l_k(-1)f^{\SC}_{\tilde{L}k} + f^{\RS}_L\right) + \sum_{\alpha=0}^i l_{\alpha}(+1)\left(f^{\SC}_{\alpha\tilde{R}} - \sum_{k=0}^N l_k(+1)f^{\SC}_{\tilde{R}k} + f^{\RS}_R\right).
\end{lequation}
Now we set~$i=N$ and apply different simplifications to prove that~\cref{eq:fc_definition_L} is equivalent to $\bar{f}_{N+1}=f^{\RS}_R$\three{\st{ when~$i=N$}}. \three{In this case}, the third and fourth terms of the right-hand side can be further simplified by considering that, for the Lagrange interpolating polynomials,~$\sum_{i=0}^N l_i(x) = 1$. Substituting this allows us to ``exchange''~$f^{\RS}_L$ in the first term by~$f^{\RS}_R$,
\begin{lequation*}
    \bar{f}_{N+1} = f^{\RS}_R + \sum_{\alpha=0}^N\sum_{k=0}^N \one{S_{\alpha k}}f^{\SC}_{\alpha k} - \sum_{k=0}^N \left[l_k(+1)f^{\SC}_{\tilde{R}k} - l_k(-1)f^{\SC}_{\tilde{L}k}\right] + \sum_{\alpha=0}^N \left[l_{\alpha}(+1)f^{\SC}_{\alpha\tilde{R}} - l_{\alpha}(-1)f^{\SC}_{\alpha\tilde{L}}\right].
\end{lequation*}
Written in this form we can now apply the symmetry property of the entropy-conservative numerical flux, $f^{\SC}_{ij}=f^{\SC}_{ji}$, and cancel the last two terms,
\begin{lequation*}
    \bar{f}_{N+1} = f^{\RS}_R + \sum_{\alpha=0}^N\sum_{k=0}^N \one{S_{\alpha k}}f^{\SC}_{\alpha k}.
\end{lequation*}
Finally, we remark that $\one{S_{ij}}=-\one{S_{ji}}$ and thus, the product of the last term is the Hadamard product of the skew-symmetric matrix~$\one{\mat{S}}$ with the symmetric matrix~$f^{\SC}_{ij}$. This is, in fact, another skew-symmetric matrix and the last term is also zero, $\bar{f}_{N+1} = f^{\RS}_R$.
\end{proof}

Note that this proof is also applicable to the \textit{standard} Gauss DGSEM scheme since the expression for the staggered fluxes is similar, simply lacking the terms containing entropy-projected quantities.

%---------------------------------------------------------------------------------------------------
\section{Subcell limiting} \label{sec:subcell}
%---------------------------------------------------------------------------------------------------

The existence of an underlying flux-differencing formula for the Gauss-DGSEM enables the use of the subcell limiting strategies presented by Rueda-Ramírez et al.~\cite{RUEDARAMIREZ2022subcell} to improve the robustness of the method, e.g., in the presence of shocks.
In particular, we propose a hybrid scheme obtained as a convex combination of the high-order Gauss-DGSEM with a first-order FV method,
\begin{lequation} \label{eq:hybrid-scheme}
    m_i \two{u_{t,i}} = \hat{f}_i - \hat{f}_{i+1},
\end{lequation}
with~$m_i$ the entries of the diagonal mass matrix. The individual fluxes are given as
\begin{lequation} \label{eq:convexComb}
    \hat{f}_i = \alpha_i \bar{f}_i^{\mathrm{FV}} + (1 - \alpha_i) \bar{f}_i, 
\end{lequation}
where~$\tilde{f}_i^{\mathrm{FV}}$ is a robust first-order approximation of the flux and $\alpha_i$ is a so-called blending coefficient, which is selected such that the resulting scheme exhibits some desired properties, e.g., positivity, non-oscillatory behavior, etc.

\two{
When using the DGSEM on Gauss-Lobatto nodes, it is possible to combine the high-order and low-order methods at the element level, i.e., with $\alpha_i$ a piece-wise constant function that is different at each element.
This limiting technique, known as element-wise blending, does not require a flux-differencing formula and has been shown to retain conservation for any choice of  $\alpha$ \cite{Hennemann2021,rueda2021entropy}.
The proofs in \cite{Hennemann2021,rueda2021entropy} rely on the fact that the surface fluxes of the low-order FV method and high-order LGL-DGSEM are equal at the boundary nodes, i.e. $\bar{f}_i^{\mathrm{FV}} = \bar{f}_i$ for $i\in \{ 0, N+1 \}$.

When using the DGSEM on Gauss nodes, we will always require the flux-differencing formula and subcell-wise blending \eqref{eq:convexComb} to ensure conservation.
Since the surface fluxes of the low-order FV method and high-order LGL-DGSEM are in general not equal at the boundary nodes, i.e. $\bar{f}_i^{\mathrm{FV}} \ne \bar{f}_i$ for $i\in \{ 0, N+1 \}$, we lose conservation properties if we select the blending coefficient of the boundary fluxes independently at each element.
As a result, we will always need to use the same blending coefficient on both sides of an inter-element interface.
}

\two{
\st{The main difference of [7] with the hybrid LGL-DGSEM/FV methods in~[10, 11] is that the Gauss-DGSEM does not share the same boundary fluxes with the first-order FV scheme. As a result, the inter-element fluxes have to be obtained as a convex combination of low- and high-order fluxes to obtain a robust scheme.}}

%---------------------------------------------------------------------------------------------------
\section{Results}
%---------------------------------------------------------------------------------------------------

This section includes some numerical results that confirm the theoretical developments of \cref{sec:theory,sec:subcell}. In both cases we solve the Euler equations using the \one{local Lax-Friedrichs numerical flux}. In \cref{sec:convergence}, the non-dissipative term uses the \one{two-point flux from Chandrashekar~\cite{Chandrashekar2013}}, whereas we employ a simple average in \cref{sec:blast}.

\subsection{Convergence} \label{sec:convergence}
%---------------------------------------------------------------------------------------------------

\begin{figure}
\centering
\noindent\begin{subfigure}{\textwidth}
    \centering
    \two{\begin{tabular}{@{} l *{5}{c} @{}}
        \toprule
        & $h = 1/2$ & $h = 1/4$ & $h = 1/8$ & $h = 1/16$ & $h = 1/32$ \\
        \midrule
        $N = 1$ & $0.20862 \cdot 10^{-13}$ & $0.14692 \cdot 10^{-13}$ & $0.40109 \cdot 10^{-13}$ & $0.10079 \cdot 10^{-12}$ & $0.20995 \cdot 10^{-12}$ \\
        $N = 2$ & $0.18708 \cdot 10^{-13}$ & $0.37276 \cdot 10^{-13}$ & $0.74427 \cdot 10^{-13}$ & $0.18228 \cdot 10^{-12}$ & $0.26580 \cdot 10^{-12}$ \\
        $N = 3$ & $0.33172 \cdot 10^{-13}$ & $0.60155 \cdot 10^{-13}$ & $0.10953 \cdot 10^{-12}$ & $0.24713 \cdot 10^{-12}$ & $0.51435 \cdot 10^{-12}$ \\
        $N = 4$ & $0.34284 \cdot 10^{-13}$ & $0.85036 \cdot 10^{-13}$ & $0.21604 \cdot 10^{-12}$ & $0.34727 \cdot 10^{-12}$ & $0.62364 \cdot 10^{-12}$ \\
        $N = 5$ & $0.47068 \cdot 10^{-13}$ & $0.11496 \cdot 10^{-12}$ & $0.20316 \cdot 10^{-12}$ & $0.41188 \cdot 10^{-12}$ & $0.73492 \cdot 10^{-12}$ \\
        \bottomrule
    \end{tabular}}
    \caption{\two{$L_2$ difference.}}
\end{subfigure}
\noindent\begin{subfigure}{0.4\textwidth}
    \centering
    \includegraphics[width=\linewidth]{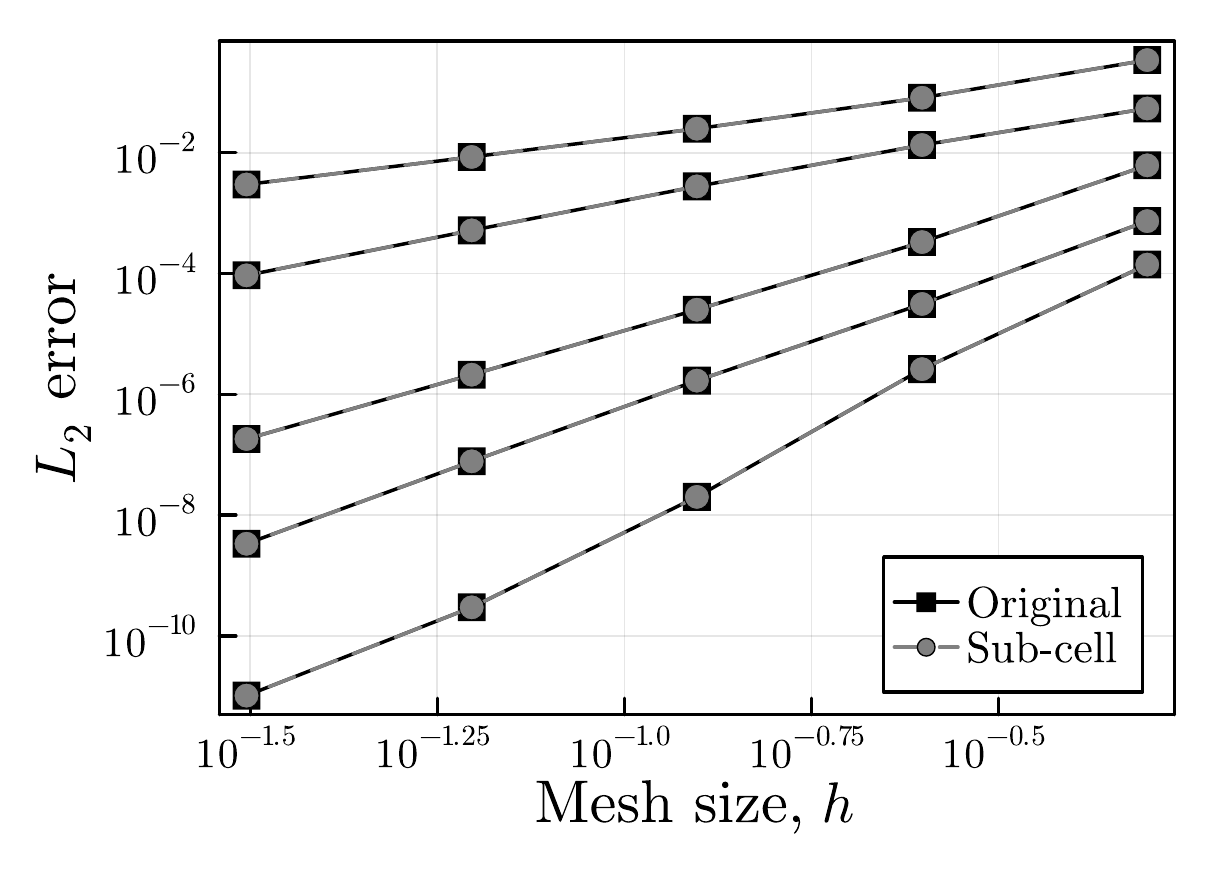}
    \caption{\two{Convergence rate.}}
\end{subfigure}
\caption{Comparison between the original formulation of J. Chan and the subcell approach with~$N = 1, \dots, 5$ (higher orders have lower errors). The two formulations overlap and cannot be distinguished in the graph.}
\label{fig:comparison}
\end{figure}

We test the numerical accuracy of the subcell approach by comparing it against the\one{\st{ split-form}} formulation obtained by Chan for Gauss nodes~\cite{chan2019efficient}. Considering the solution to the one-dimensional Euler equations,
\begin{lequation*}
    \rho = 2 + \sin \pi (x - t), \quad u = 1, \quad p = 1,
\end{lequation*}
we integrate them until $t=0.7$ with a $CFL=0.125$ and using the 5-stage, 4th-order Runge-Kutta algorithm described by Carpenter and Kennedy~\cite{Carpenter1994}. The \two{results of the comparison between the new approach and the baseline from J. Chan} are shown in~\cref{fig:comparison}.

\subsection{Sedov blast} \label{sec:blast}
%---------------------------------------------------------------------------------------------------

To illustrate the shock-capturing capacity of the hybrid DGSEM/FV method, we simulate a Sedov blast problem describing the evolution of a blast wave expanding from an initial concentration of density and pressure.
For the initial condition, we assume a gas in rest, $v_1 (t=0) = v_2 (t=0) =0$, with a Gaussian distribution of density and pressure,
\begin{lequation*}
    \rho (t=0) = \rho_0 + G(r; \sigma_{\rho}), \quad
    p (t=0) = p_0 + (\gamma - 1)G(r; \sigma_p), \quad
    G(r; \sigma) = \frac{1}{4 \pi \sigma^2} \exp \left( -\frac{1}{2} \frac{r^2}{\sigma^2} \right), \quad
    r^2 = x^2 + y^2,
\end{lequation*}
where we choose $\sigma_{\rho}=0.25$ and $\sigma_p=0.15$. Furthermore, the ambient density is set to $\rho_0=1$ and the ambient pressure to $p_0=10^{-1}$. We complement the simulation domain, $\Omega=[-1,1]^2$, with periodic boundary conditions, and tessellate it using $K=64^2$ quadrilateral elements. We represent the solution with polynomials of degree $N=3$ and run the simulation until $t=1$.

To avoid non-physical oscillations in the vicinity of shocks, we use \one{a subcell-wise limiting strategy \eqref{eq:convexComb} to impose a non-oscillatory behavior on the density of each node $i$,
\begin{lequation} \label{eq:rho_condition}
    \min_{j \in \mathcal{N}(i)} \overline{\rho}_{ij} \le \rho_i \le \min_{j \in \mathcal{N}(i)} \overline{\rho}_{ij}, \quad
    \overline{u}_{ij} = \left(\overline{\rho}_{ij},\overline{\rho \vec{v}}_{ij},\overline{\rho E}_{ij} \right)^T \coloneqq \frac{1}{2} (u_i + u_j) + \frac{\hat{\vec{n}}_{ij}}{2 \lambda^{\max}_{ij} } \cdot (f_j - f_i),
\end{lequation}
where $\mathcal{N}(i)$ is the so-called low-order stencil of node $i$, a set containing all the neighboring nodes to $i$, $\overline{u}_{ij}$ is the so-called \textit{bar state}, $\lambda^{\max}_{ij}$ is an estimate of the maximum wave speed between nodes $i$ and $j$, and $\hat{\vec{n}}_{ij}$ is the normal vector (normalized metric term) at the interface between nodes $i$ and $j$.
We refer the reader to \cite{RUEDARAMIREZ2022subcell} for details on how to compute $\alpha$ with an algebraic flux correction scheme, such that \eqref{eq:rho_condition} is guaranteed.

Using the procedure described in \cite{RUEDARAMIREZ2022subcell}, we obtain a provisional blending coefficient $\tilde \alpha_i$ at each node of our domain.
The subcell-wise limiting strategy \eqref{eq:convexComb} is then applied taking the maximum blending coefficient on both sides of every interface, $\alpha_i = \max (\tilde \alpha_i, \tilde \alpha_{i+1})$.
The same procedure is applied at element interfaces when using Gauss nodes.
As explained above, the surface terms of DG and FV are equal when using Gauss-Lobatto nodes.
Hence, it is not necessary to blend the fluxes at the interfaces.
}

\begin{figure}[h!]
    \centering
    \begin{overpic}[trim=164 474 882 29,clip,width=0.2\linewidth]{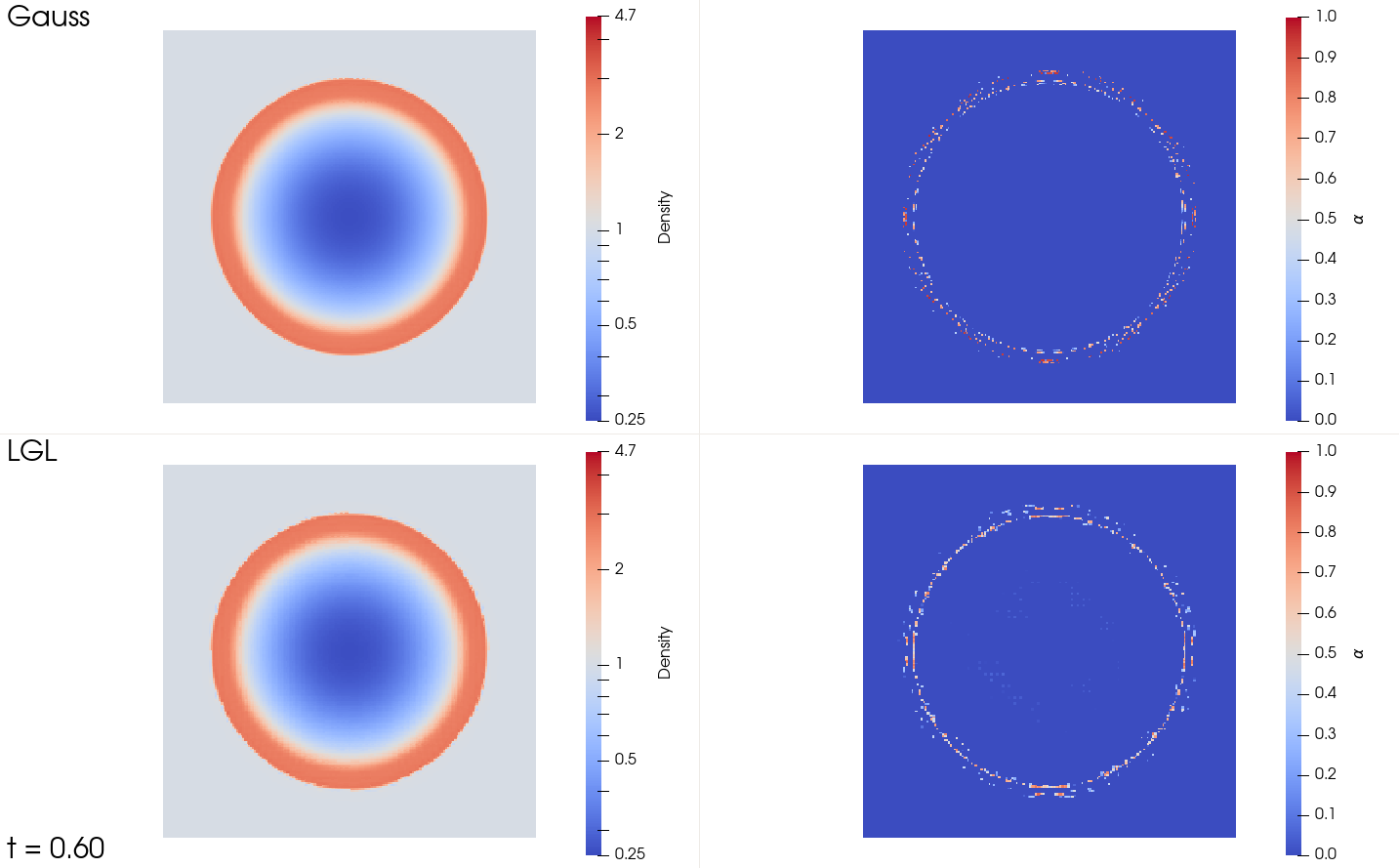}
      \put(1,1) {\small Gauss}
    \end{overpic}
    \begin{overpic}[trim=164 29 882 474,clip,width=0.20\linewidth]{figures/SC_Gauss_vs_LGL.0090.png}
      \put(1,1) {\small LGL}
    \end{overpic}
    \includegraphics[trim=595 452 740 0,clip,height=0.20\linewidth]{figures/SC_Gauss_vs_LGL.0090.png}
    \hfill
    \begin{overpic}[trim=882 474 164 29,clip,width=0.2\linewidth]{figures/SC_Gauss_vs_LGL.0090.png}
      \put(1,1) {\small\color{white} Gauss}
    \end{overpic}
    \begin{overpic}[trim=882 29 164 474,clip,width=0.2\linewidth]{figures/SC_Gauss_vs_LGL.0090.png}
      \put(1,1) {\small\color{white} LGL}
    \end{overpic}
    \includegraphics[trim=1311 0 0 452,clip,height=0.20\linewidth]{figures/SC_Gauss_vs_LGL.0090.png}
    
    \caption{\one{Density and blending coefficient,~$\alpha$, for the Sedov blast problem obtained with the hybrid DGSEM/FV method at $t=0.60$ using Gauss and Gauss-Lobatto nodes.}}
    \label{fig:blast}
\end{figure}

\Cref{fig:blast} illustrates the density and blending coefficient in the domain obtained with the hybrid DGSEM/FV method using Gauss and Gauss-Lobatto nodes at time $t=0.6$.
Both schemes capture the expanding shock correctly while using the high-order DG method in most of the domain.

\one{
We compute the mean blending coefficient at a particular time as 
\begin{lequation*}
    \bar \alpha (t) = \left( \frac{1}{V} \sum_{e=1}^K \sum_{i,j=0}^N J_{ij} \omega_{ij} \hat \alpha^e_{ij}(t) \right),
\end{lequation*}
where $e \in [1,K]$ denotes the element index, $K$ is the number of elements of the domain, $i,j \in [0,N]$ are the node indices, $N$ is the polynomial degree, $\hat \alpha^e_{ij}(t)$ is the provisional (nodal)  blending coefficient of node $ij$ of element $e$ at time $t$, and $V$ is the total area of the domain.

\Cref{fig:evol} shows the evolution of the mean blending coefficient and number of time steps taken for the simulation of the blast wave with CFL=0.9.
For this particular example and the non-oscillatory condition on density \eqref{eq:rho_condition}, the ES Gauss scheme requires more limiting than the ES LGL scheme for the same number of degrees of freedom, polynomial degree, and CFL number.
However, the LGL scheme takes more time steps to finish the simulation as it requires shorter time-step sizes.

The maximum allowable time-step size that is required for condition \eqref{eq:rho_condition} in 2D reads \cite{rueda2023monolithic}
\begin{lequation*}
    \Delta t \le \min_{ij}
    \frac{J\,\omega_i\,\omega_j}{\omega_j \left(\lambda^{\max}_{(i,i-1)j}+\lambda^{\max}_{(i,i+1)j}\right) + \omega_i \left(\lambda^{\max}_{i(j,j-1)}+\lambda^{\max}_{i(j,j+1)}\right)}.
\end{lequation*}
The larger quadrature weights of the Gauss-DGSEM at the boundary nodes of the elements result in a less restrictive CFL condition.
As a result, the DGSEM simulation with Gauss nodes is not only more accurate, but also computationally more efficient than the LGL simulation for the same number of degrees of freedom.

\begin{figure}
    \centering
    \begin{subfigure}{0.4\linewidth}
    \includegraphics[width=\textwidth]{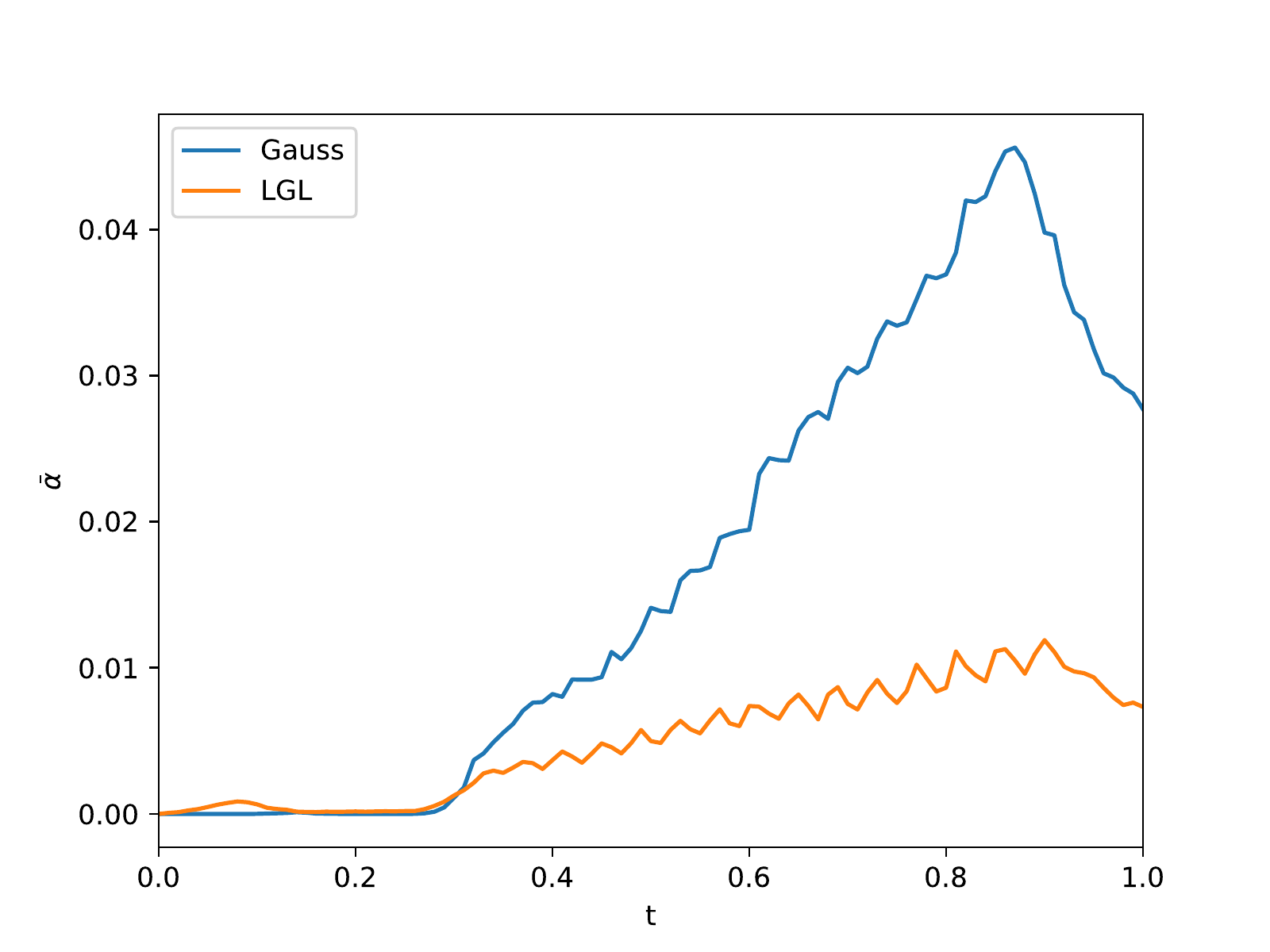}
    \caption{\one{Mean blending coefficient}}
    \end{subfigure}
    \begin{subfigure}{0.4\linewidth}
    \includegraphics[width=\textwidth]{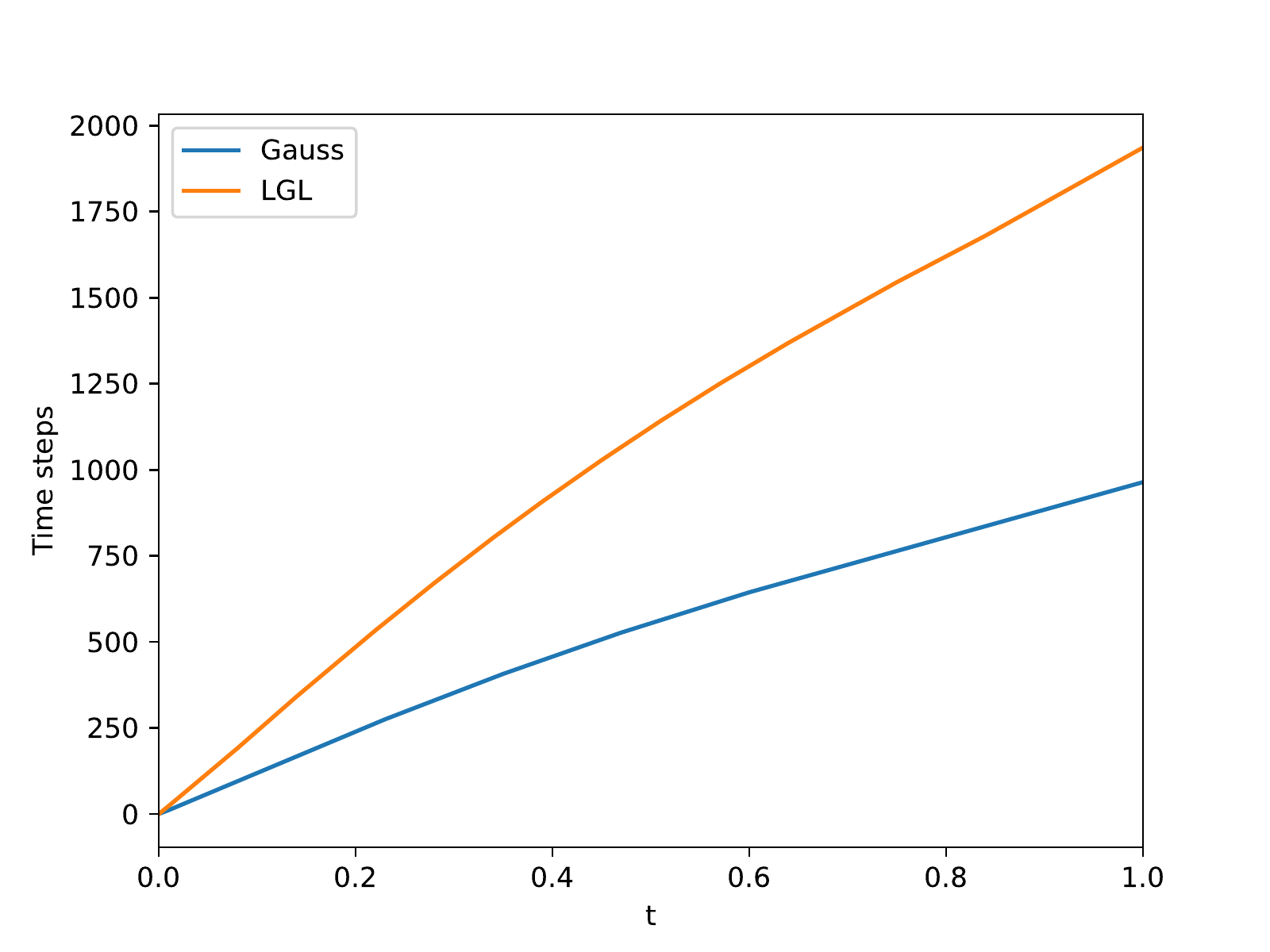}
    \caption{\one{Number of time steps}}
    \end{subfigure}
    \caption{\one{Evolution of the mean blending coefficient and number of time steps taken for the simulation of the blast wave with CFL=0.9.}}
    \label{fig:evol}
\end{figure}
}
%---------------------------------------------------------------------------------------------------
\section{Conclusions}
%---------------------------------------------------------------------------------------------------

We have presented in this work a novel flux-differencing expression of the entropy-conservative formulation of J. Chan~\cite{Chan2018,chan2019efficient} in terms of staggered fluxes. This approach can be applied explicitly to the DGSEM framework with Gauss nodes. We have also proved that this telescopic formulation of the derivative operators exists\one{\st{ for split-forms}} with and without entropy-projected variables.

When applied to a simple test case, we have shown that the convergence properties of the flux-differencing method match those of the entropy-conservative \one{approach} used as the baseline. This allowed us to apply some subcell limiting strategies already developed for the DGSEM with Gauss-Lobatto nodes, although more work is needed in this aspect.

%---------------------------------------------------------------------------------------------------
\section{Acknowledgements}
%---------------------------------------------------------------------------------------------------
Andr\'es Mateo has received funding from Universidad Polit\'ecnica de Madrid under the Programa Propio PhD programme. Andrés M. Rueda-Ramírez acknowledges funding through the Klaus-Tschira Stiftung via the project ``HiFiLab''. \one{Gonzalo Rubio and Eusebio Valero acknowledge the funding from the European Union’s Horizon 2020 research and innovation programme under the Marie Skłodowska Curie grant agreements No 955923-SSECOID and 101019137-FLOWCID}. We furthermore thank the Regional Computing Center of the University of Cologne (RRZK) for providing computing time on the High Performance Computing (HPC) system ODIN as well as support. Finally, all authors gratefully acknowledge Universidad Politécnica de Madrid (www.upm.es) for providing computing resources on the Magerit Supercomputer.

\appendix

\three{

%---------------------------------------------------------------------------------------------------
\section{Extension to higher dimensions and curvilinear grids}
%---------------------------------------------------------------------------------------------------

The flux-differencing formula \eqref{eq:split_chan_fv} extends to multiple space dimensions on curvilinear grids. 
We can write the Gauss-DGSEM in two dimensions as
\begin{lequation} \label{eq:fv_2d}
    J_{ij} \frac{\partial u_{ij}}{\partial t} = \frac{1}{\omega_i} \left(\bar{f}_{(i-1,i)j} - \bar{f}_{(i,i+1)j}\right) + \frac{1}{\omega_j} \left(\bar{f}_{i(j-1,j)} - \bar{f}_{i(j,j+1)} \right),
\end{lequation}
where we introduce a new notation: $\bar{f}_{(i-1,i)j}$ denotes the high-order telescoping flux between nodes $(i-1,j)$ and $(i,j)$, and $\bar{f}_{i(j-1,j)}$ denotes the high-order telescoping flux between nodes $(i,j-1)$ and $(i,j)$.

The telescoping fluxes are uniquely defined between two neighboring nodes. For instance, in the first coordinate direction they read,
\begin{lequation}\label{eq:teleflux_2d}
\begin{alignedat}{3}
    \bar{f}_{(-1,0)j} &= f^{\RS}_{Lj}, \\
    \bar{f}_{(i,i+1)j} &= \bar{f}_{(i-1,i)j} + \sum_{k=0}^N S_{ik}f^{1\SC}_{(i,k)j} &&- l_i(-1)\left[f^{1\SC}_{(i,\tilde{L})j} - \sum_{k=0}^N l_k(-1)f^{1\SC}_{(\tilde{L},k)j} + f^{\RS}_{Lj} \right] \\
        &&&+ l_i(+1)\left[f^{1\SC}_{(i,\tilde{R})j} - \sum_{k=0}^N l_k(+1)f^{1\SC}_{(\tilde{R},k)j} + f^{\RS}_{Rj} \right], \quad i=0, \ldots, N,
\end{alignedat}
\end{lequation}
where the surface numerical fluxes depend on the inner entropy-projected solution, the outer solution, and the normal vector at the boundary of the element,
\begin{lequation*}
\begin{aligned}
    f^{\RS}_{Lj} &= f^{\RS}_{Lj} \left( u(v_{Lj}),u^+_{Lj},\vec{n}_{Lj} \right), \\
    f^{\RS}_{Rj} &= f^{\RS}_{Rj} \left( u(v_{Rj}),u^+_{Rj},\vec{n}_{Rj} \right).
\end{aligned}
\end{lequation*}
The volume numerical fluxes are entropy-conservative two point fluxes with metric dealiasing (see, e.g., \cite{Rueda2023}),
\begin{lequation*}
    f^{1\SC}_{(i,k)j} = f^{\SC} \left( u_{ij}, u_{kj} \right) \cdot \avg{J\vec{a}^1}_{(i,k)j},
\end{lequation*}
and we use the so-called contravariant metric vector, $\vec{a}^1 \coloneqq \nabla \xi$, which relates the physical-frame coordinates $(x,y)$ with the reference-frame coordinates $(\xi, \eta)$, and the average operator,
\begin{lequation*}
    \avg{a}_{(i,j)} \coloneqq \frac{1}{2} (a_i + a_j).
\end{lequation*}
The telescoping fluxes are defined analogously in the other coordinate directions.

\subsection{Subcell limiting}

Due to the existence of a flux-differencing formula in multiple space dimensions for the Gauss-DGSEM \eqref{eq:fv_2d}, it is possible to apply subcell limiting strategies. However, we need the subcell metric terms to compute the low-order fluxes. Following the strategy presented by~\citet{Hennemann2021}, we replace a constant state in~\eqref{eq:teleflux_2d} to obtain the subcell metric terms:
\begin{lequation*}
\begin{alignedat}{3}
    \vec{n}_{(-1,0)j} &= \vec{n}_{Lj}, \\
    \vec{n}_{(i,i+1)j} &= \vec{n}_{(i-1,i)j} + \sum_{k=0}^N S_{ik} \avg{J\vec{a}^1}_{(i,k)j} &&- l_i(-1)\left[\avg{J\vec{a}^1}_{(i,\tilde{L})j} - \sum_{k=0}^N l_k(-1)\avg{J\vec{a}^1}_{(\tilde{L},k)j} + \vec{n}_{Lj} \right] \\
    &&&+ l_i(+1)\left[\avg{J\vec{a}^1}_{(i,\tilde{R})j} - \sum_{k=0}^N l_k(+1)\avg{J\vec{a}^1}_{(\tilde{R},k)j} + \vec{n}_{Rj} \right]&, \quad i=0, \ldots, N,
\end{alignedat}
\end{lequation*}
which simplifies to
\begin{lequation*}
\begin{alignedat}{3}
    \vec{n}_{(-1,0)j} &= \vec{n}_{Lj}, \\
    \vec{n}_{(i,i+1)j} &= \vec{n}_{(i-1,i)j} + \sum_{k=0}^N S_{ik} \avg{J\vec{a}^1}_{(i,k)j} &&- l_i(-1)\left[ \frac{1}{2}\left(J\vec{a}^1\right)_{ij} - \frac{1}{2} \sum_{k=0}^N l_k(-1) \left(J\vec{a}^1\right)_{kj} + \vec{n}_{Lj} \right] \\
    &&&+ l_i(+1)\left[\frac{1}{2} \left(J\vec{a}^1\right)_{ij} - \frac{1}{2} \sum_{k=0}^N l_k(+1)\left(J\vec{a}^1\right)_{kj} + \vec{n}_{Rj} \right]&, \quad i=0, \ldots, N.
\end{alignedat}
\end{lequation*}
}

%---------------------------------------------------------------------------------------------------
\bibliographystyle{model1-num-names}
\bibliography{library}
%---------------------------------------------------------------------------------------------------

%---------------------------------------------------------------------------------------------------
\end{document}